\documentclass[reqno,12pt]{amsart}
\usepackage{graphicx,epsfig}

\usepackage{amssymb,amsthm,amsmath}
\usepackage{latexsym}

\newtheorem{definition}{Definition}
\newtheorem{thm}{Theorem}

\newtheorem{rem}[thm]{Remark}

\newcommand\cal{\mathcal}

\def\Z{\hbox{\font\dubl=msbm10 scaled 1200 {\dubl Z}}}
\def\R{\hbox{\font\dubl=msbm10 scaled 1200 {\dubl R}}}
\def\Q{\hbox{\font\dubl=msbm10 scaled 1200 {\dubl Q}}}
\def\K{\hbox{\font\dubl=msbm10 scaled 1200 {\dubl K}}}

\title[Right triangles with algebraic sides and elliptic curves ]
{Right triangles with algebraic sides and elliptic curves over number fields}

\author[E. Girondo et al.]{E. Girondo, G. Gonz\'alez-Diez, E. Gonz\'alez-Jim\'enez,\\ R. Steuding, J. Steuding}

\subjclass[2000]{Primary 11G05; Secondary 11A99}
\date{\today}

\begin{document}

\maketitle

\begin{abstract}
Given any positive integer $n$, we prove the existence of infinitely many right triangles with area $n$ and side lengths in certain number fields. This generalizes the famous congruent number problem. The proof allows the explicit construction of these triangles; for this purpose we find for any positive integer $n$ an explicit cubic number field $\Q(\lambda)$ (depending on $n$) and an explicit point $P_\lambda$ of infinite order in the Mordell-Weil group of the elliptic curve $Y^2=X^3-n^2X$ over $\Q(\lambda)$.
\end{abstract}

\section{Congruent numbers over the rationals}
A positive integer $n$ is called a congruent number if there exists a right triangle with rational sides and area equal to $n$, i.e., there exist $a,b,c\in\Q^*$ with
\begin{equation}\label{one}
a^2+b^2=c^2\qquad\mbox{and}\qquad {\textstyle{1\over 2}}ab=n.
\end{equation}
It is easy to decide whether there is a right triangle of given area and integral sides (thanks to Euclid's characterization of the Pythagorean triples). The case of rational sides, known as the congruent number problem, is not completely understood. Fibonacci showed that $5$ is a congruent number (since one may take $a={3\over 2}$, $b={20\over 3}$ and $c={41\over 6}$). Fermat found that $1$, $2$ and $3$ are not congruent numbers. Hence, there is no perfect square amongst the congruent numbers since otherwise the corresponding rational triangle would be similar to one with area equal to $1$. By the same reason one can restrict to square-free positive integers, a condition we will assume in the sequel.

\par

There is a fruitful translation of the congruent number problem into the language of elliptic curves. Suppose $n$ is a congruent number. It follows from (\ref{one}) that there exist three squares in arithmetic progression of distance $n$, namely $x-n,x,x+n$, where $x=c^2/4$. Therefore their product $(x-n)x(x+n)$ is again a rational square. In other words, a right triangle of area $n$ and rational sides $a$, $b$, $c$ corresponds to the rational point $(c^2/4,c(a^2-b^2)/8)$ on the elliptic curve
$$
E_n\,:\,Y^2=X^3-n^2X.
$$
Conversely, given a rational point $(x,y)$ on $E_n$ such that $y\neq 0$, one may define
\begin{equation}\label{two}
a=\left\vert {y\over x}\right\vert, \quad b=2n\left\vert{x\over y}\right\vert, \quad c={x^2+n^2\over \vert y\vert}
\end{equation}
to obtain a right triangle with rational sides $a$, $b$, $c$ and area $n$.

Therefore, $n$ is a congruent number if and only if the elliptic curve $E_n$ has a rational point with non-vanishing $y$-coordinate. The correspondence between rational points on $E_n$ and right triangles with rational sides and area $n$ is not bijective. For instance, solving $x$ and $y$ with respect to $a,b$ and $c$ in equation (\ref{two}) gives the two points
\begin{equation}\label{three}
x={1\over 2}a(a\pm c),\qquad y=ax\,.
\end{equation}

The points $(x,y)$ with $y\neq 0$ have infinite order in the Mordell-Weil group $E_n(\Q)$, since it is known (see \cite{silverman}) that the torsion subgroup of $E_n(\Q)$ consists only of points of order $2$, that are $(0,0), (\pm n,0)$, and the point at infinity. From this one can deduce the fact, already known to Fermat, that there are in fact infinitely many right triangles with rational sides $a,b,c$ verifying equation (\ref{one}) for a given congruent number $n$, since scalar multiplication of the points corresponding to (\ref{three}) on $E_n$ yields new right triangles of area $n$.

\

Tunnell \cite{tunnell} found an intruiging approach towards the congruent number problem. He showed that if $n$ is an odd square-free congruent number then
\begin{eqnarray*}
\lefteqn{\sharp\{x,y,z\in\Z\,:\,2x^2+y^2+8z^2=n\}}\\
&=&2\sharp\{x,y,z\in\Z\,:\,2x^2+y^2+32z^2=n\},
\end{eqnarray*}
and that the converse implication is also true provided the yet unsolved Birch \& Swinnerton-Dyer conjecture holds (i.e., the rank of an elliptic curve $E$ is positive if and only if the associated $L$-function vanishes at the central point: $L(E,1)=0$). A similar criterion exists when $n$ is even. This conjectural equivalent statement for the congruent number problem in terms of the number of representations of $n$ by certain ternary quadratic forms allows to determine by computation whether a given integer $n$ is congruent or not.

\section{Congruent numbers over number fields}

A natural generalization of the congruent number problem is to allow the sides of the triangles to belong to an algebraic number field. The idea to study the congruent number problem over algebraic extensions dates back at least to Tada \cite{tada} who considered real quadratic fields.

\begin{definition}
 Let $\K$ be a number field and let $n$ be a positive integer. We say that $n$ is a $\K$-congruent number if there exist $a,b,c\in\K$ such that (\ref{one}) holds.
\end{definition}

\noindent Note that $\Q$-congruent numbers are nothing but the usual congruent numbers. When $\K$ is a subfield of $\R$ the geometric meaning still holds, since $n$ being $\K$-congruent implies the existence of a right triangle with real sides in $\K$ and area equal to $n$.

It is easy to see that any positive integer $n$ is a congruent number for some quadratic extension of $\Q$. For instance, $a=b=\sqrt{2}$ yields a right triangle with area equal to $n=1$. However, equation (\ref{three}) leads to the points $(1\pm\sqrt{2},\sqrt{2}\pm 2)$ on the curve $E_1$ which are all torsion points over $\Q(\sqrt{2})$. Therefore we do not get infinitely many different right triangles from these points. This example motivates the following definition:

\begin{definition}
 We say that a $\K$-congruent number $n$ is \emph{properly} $\K$-congruent if (\ref{one}) has infinitely many solutions $a,b,c \in \K$.
\end{definition}

\noindent This definition is not interesting for $\K= \Q$, since all $\Q$-congruent numbers are properly $\Q$-congruent. But the example above shows that things may be different over number fields.

\subsection{Congruent numbers over quadratic fields}

\begin{thm}
Every positive integer $n$ is properly $\K$-congruent over some real quadratic field $\K$.
\end{thm}

\par\medskip

\noindent {\bf Proof}. Consider the system
$$
\left\{  \begin{array}{l}
a^2+b^2= c^2\,, \\
ab =2n\,.
\end{array} \right.
$$
Substituting $a$ by $2n/b$ in the first equation yields
$$
4n^2 + b^4 = c^2 b^2.
$$
Therefore
$$
c = \frac{\sqrt{4n^2+ b^4}}{b}.
$$
It follows that $n$ is congruent over $\mathbb{Q}(\sqrt{m})$, where  $m=4 n^2+b^4$ (and $b$ is any chosen rational point).

\

Given such a triple $(a,b,c) \in \mathbb{Q}(\sqrt{m})^3 $, the corresponding $\mathbb{Q}(\sqrt{m})$-rational point $P$ on the curve $E_n$ is given by (\ref{three}). It is known (see \cite{tada,kwon}) that the torsion subgroup of $E_n(\mathbb{Q}(\sqrt{m}))$ reduces to the 2-torsion, namely
$\{ \infty, (0,0), (\pm n, 0) \}$, except for the cases $(n,m)=(1,2)$ or $(n,m)=(2,2)$.

\

We can choose $b \in  \mathbb{Q}$ such that $m \neq 2$, and therefore $P$ is a non-torsion point. This finishes the proof.
\par\medskip

\begin{rem} \rm
Tada \cite{tada} showed that if $n$ is not $\mathbb{Q}$-congruent, then it is congruent over $\mathbb{Q}(\sqrt{m})$ if and only if $n m$ is $\mathbb{Q}$-congruent. This fact has applications of the following kind.

Since it is known that 1 and 11 are not $\mathbb{Q}$-congruent numbers, it follows that 1 is not $\mathbb{Q}(\sqrt{11})$-congruent. According to the proof above, $4 n^2 + b^4=4 + b^4$ cannot be of the form $11 \cdot d^2$ with $d \in \mathbb{Q}$ for any choice of $b \in \mathbb{Q}$. That is, the genus 1 curve
$$
11 y^2 = x^4 + 4
$$
is not an elliptic curve over $\mathbb{Q}$, i.e. it has no nonsingular $\mathbb{Q}$-rational point.

\

The same argument applied to a general pair $n,m\in\Z^*$ leads to the following generalization. Denote by $C_{n,m}$ the genus 1 curve $ C_{n,m}: my^2 = x^4 + 4 n^2$. Then this curve has a $\Q$-rational point if and only if $n$ is $\Q(\sqrt{m})$-congruent, or equivalently, if and only if $n m$ is  $\Q$-congruent.

Note that the curve $ C_{n,m}$ is the twist by $m$ of the curve $ C_{n,1}: y^2 = x^4 + 4 n^2$, that is another parametrization of the congruent number problem.

\end{rem}

\subsection{Congruent numbers over cubic fields}

Next we shall construct an explicit non-torsion point on $E_n$ over certain cubic fields (depending on $n$).

\begin{thm}\label{3}
Every positive integer $n$ is properly congruent over some real cubic field. More precisely, for a positive integer $n$ let 
$\lambda=\lambda(n)$ be the unique real solution of the cubic equation 
$$
32\lambda^3 - 32\lambda^2 +8\lambda + n^2=0.
$$
Then the point $P_\lambda=(x_\lambda,y_\lambda)$ with coordinates given by (\ref{coord}) and (\ref{coord1}) below is of infinite order in the Mordell-Weil group of the elliptic curve $Y^2=X^3-n^2X$ over $\Q(\lambda)$.
\end{thm}

\noindent {\bf Proof}. We use an idea of Chahal \cite{chahal}. Starting from an old identity of Desboves \cite{desboves},
\begin{eqnarray*}
\lefteqn{({\cal Y}^2+2{\cal X}{\cal Y}-{\cal X}^2)^4+(2{\cal X}^3{\cal Y}+{\cal X}^2{\cal Y}^2)(2{\cal X}+2{\cal Y})^4}\\
&=&({\cal X}^4+{\cal Y}^4+10{\cal X}^2{\cal Y}^2+4{\cal X}{\cal Y}^3+12{\cal X}^3{\cal Y})^2,
\end{eqnarray*}
the substitution ${\cal X}=1-2\lambda, {\cal Y}=4\lambda$ yields
\begin{eqnarray*}
\lefteqn{(1-12\lambda+4\lambda^2)^4+8\lambda(2\lambda-1)^2(2(1+2\lambda))^4}\\
&=&(1+40\lambda-104\lambda^2+160\lambda^3+16\lambda^4)^2.
\end{eqnarray*}
This gives the point $(x,y)$ with coordinates
\begin{eqnarray}
x&=&x(\lambda)={(1-12\lambda+4\lambda^2)^2\over 4(1+2\lambda)^2},\label{tzt1}\\
y&=&y(\lambda)={(1-12\lambda+4\lambda^2)(1+40\lambda-104\lambda^2+160\lambda^3+16\lambda^4)
\over 8(1+2\lambda)^3}\label{tzt2}
\end{eqnarray}
on the elliptic curve $y^2=x^3+dx$, where $d=d(\lambda)=8\lambda(2\lambda-1)^2$. Now let $n$ be a positive integer and consider the equation $-n^2=8\lambda(2\lambda-1)^2$. It is easy to see that there is a unique real solution, explicitly given by
$$
\lambda=\lambda(n)={1\over 3}+{1\over 12}\kappa+{1\over 3\kappa},
$$
where
$$
\kappa=\kappa(n)=\sqrt[3]{-8-27n^2+3\sqrt{48n^2+81n^4}}.
$$
Here we choose $\kappa$ as the unique real third root which is negative for all $n$. By computation, we obtain a point $P_\lambda=(x_\lambda,y_\lambda)$ on $E_n(\Q(\lambda))$, where the coordinates are given by
\begin{eqnarray}
x=x_\lambda&=&{1\over 4(n^2-16)^2}\{256+992n^2+65n^4+(1024-2688n^2-28n^4)\lambda\nonumber \\
&&+(1024+1920n^2+4n^4)\lambda^2\},\label{coord}\\
y=y_\lambda&=&{1\over 32(n^2-16)^3}\{-16384+72704n^2+80960n^4+2868n^6-n^8\nonumber\\
&&+(196608-462848n^2-145152n^4-1456n^6)\lambda\label{coord1}\\
&&+(196608+421888n^2+100608n^4+208n^6)\lambda^2\}.\nonumber
\end{eqnarray}
This can be checked by use of a computer algebra package or by hand computation as follows: firstly, multiply the numerator and denominator of the coordinates (\ref{tzt1}) and (\ref{tzt2}) by the conjugates of $(1+2\lambda)^2$ and $(1+2\lambda)^3$, respectively, then the coordinates are given as rational polynomials in the variable $\lambda$; secondly, since $\Q(\lambda)$ is a cubic field, these polynomials turn out to have degree at most two, which yields (\ref{coord}) and (\ref{coord1}).

It remains to show that the point $P_\lambda$ has infinite order on $E_n(\Q(\lambda))$.  Note that  the map that sends $(x,y)$ to $(-x,\sqrt{-1}y)$ is an endomorphism on the elliptic curve $E_n$, hence $E_n$ has complex multiplication by $\mathbb{Z}[\sqrt{-1}]$. There are bounds \cite{boundcm,silverberg} for the order of the subgroup of torsion points  on a elliptic curve with complex multiplication defined over a number field. Applying these results to the elliptic curve $E_n$ over the cubic number field $\Q(\lambda)$ we obtain that if $M$ is the order of a torsion point of $E_n(\Q(\lambda))$, then $\phi(M)\leq 6$, where $\phi$ is Euler's $\phi$-function. Therefore $M\in\mathcal{B}=\{1,2,3,4,5,6,7,8,9,10,12,14,18\}$. Hence it is sufficient to prove that the order of $P_\lambda$ is not in $\mathcal{B}$.  For this purpose we use the $m$-division polynomial $\Psi_m(x)$ associated with the elliptic curve $E_n$. It suffices to check that $x_\lambda$ is not a root of $\Psi_m(x)$ for $m\in\mathcal{B}$ over $\Q(\lambda)$. We have
$$
\Psi_m(x_\lambda)=P_{0,m}(n)+P_{1,m}(n)\lambda+P_{2,m}(n)\lambda^2\,,
$$
where $P_{k,m}\in \Q[n]$ has no integer roots for $k=0,1,2$ and $m\in\mathcal{B}$. Thus $\Psi_m(x_\lambda)\ne 0$ for any integer $n$ and $m\in \mathcal{B}$. This  proves that $P_{\lambda(n)}$ has infinite order and so the Mordell-Weil group of $E_n(\Q(\lambda(n)))$ has positive rank. The theorem is proved.
\par\medskip

\subsubsection{An example}

The integer $n=1$ is a congruent number over the cubic field $\Q(\alpha)$, where
$$
\alpha={1\over 3}+{1\over 12}\sqrt[3]{-35+3\sqrt{129}}+{1\over 3\sqrt[3]{-35+3\sqrt{129}}}.
$$
The following picture shows a right triangle with area 1 and sides in $\Q(\alpha)$:
\begin{figure}[!htbp]
\begin{center}
\includegraphics[width=400pt]{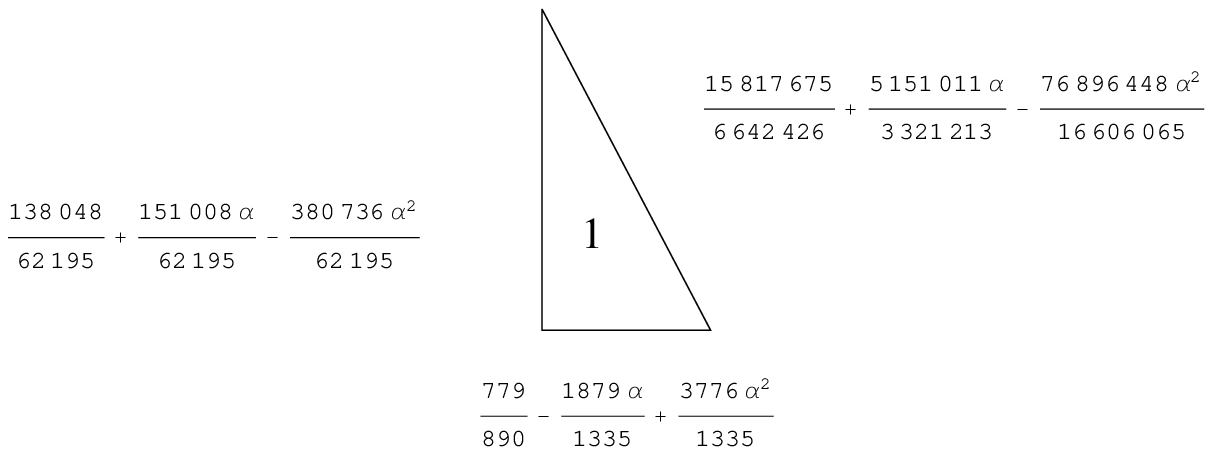}
\end{center}
\end{figure}

We conclude with an interesting but difficult question:
\begin{center}
{\it Given a number field $\K$, what are the congruent numbers over $\K$? }
\end{center}

\vspace*{.5cm}

\subsection*{Acknowledgements}
We thank J. Lehnert (Frankfurt) for useful discussions and P. M\"uller (W\"urzburg) for providing an alternative argument for Theorem 3. We also thank the referee for careful reading of the first version of this article and for pointing out some inaccuracies.
\smallskip

The fourth and the fifth author wish to thank the organizers of the 18th Czech and Slovak International Conference on Number Theory at Smolenice, August 27-31, 2007, for their kind hospitality and support.
\smallskip

Research of the third author was supported in part by grant MTM 2006-10548 (Ministerio de Educaci\'on y Ciencia, Spain) and CCG06-UAM/ESP-0477 (Universidad Aut\'onoma de Madrid - Comunidad de Madrid, Spain). Research of the rest of authors was supported in part by grant MTM 2006-01859  (Ministerio de Educaci\'on y Ciencia, Spain).

{\small

\par\bigskip

\noindent
Ernesto Girondo\\
Departamento de Matem\'aticas, CSIC\\
Madrid, Spain\\
ernesto.girondo@uam.es
\par\medskip

\noindent
Gabino Gonz\'alez-Diez, Enrique Gonz\'alez-Jim\'enez\\
Departamento de Matem\'aticas, Universidad Aut\'onoma de Madrid\\
C. Universitaria de Cantoblanco, 28\,049 Madrid, Spain\\
gabino.gonzalez@uam.es, enrique.gonzalez.jimenez@uam.es
\par\medskip

\noindent Rasa Steuding, J\"orn Steuding\\
Institut f\"ur Mathematik, Universit\"at W\"urzburg\\
Am Hubland, 97\,218 W\"urzburg, Germany\\
steuding@mathematik.uni-wuerzburg.de
}

\end{document}